\pdfoutput=1
\RequirePackage{ifpdf}
\ifpdf 
\documentclass[pdftex]{sigma}
\else
\documentclass{sigma}
\fi

\usepackage{settings}

\begin{document}

\allowdisplaybreaks

\newcommand{\arXivNumber}{2505.22607}

\renewcommand{\PaperNumber}{011}

\FirstPageHeading

\ShortArticleName{Contraction of the $\mathfrak{sl}_2$-Triple Associated to the $(k, a)$-Generalized Fourier Transform}
\ArticleName{Contraction of the $\boldsymbol{\mathfrak{sl}_2}$-Triple Associated\\ to the $\boldsymbol{(k, a)}$-Generalized Fourier Transform}

\Author{Tatsuro HIKAWA}

\AuthorNameForHeading{T.~Hikawa}

\Address{Graduate School of Mathematical Sciences, The University of Tokyo, \\
3-8-1 Komaba, Meguro-ku, Tokyo, 153-8914, Japan}
\Email{\mail{hikawa@ms.u-tokyo.ac.jp}}

\ArticleDates{Received May 29, 2025, in final form January 15, 2026; Published online February 08, 2026}

\Abstract{Ben~Sa\"{\i}d--Kobayashi--{\O}rsted introduced a family of $ \mathfrak{sl}_2 $-triples of differential-difference operators $ \mathbb{H}_{k,a} $, $ \mathbb{E}^+_{k,a} $ and $ \mathbb{E}^-_{k,a} $ on $ \mathbb{R}^N \setminus \{0\} $ indexed by a~Dunkl parameter $ k $ and a~deformation parameter $ a \neq 0 $. In the present paper, we study the behavior as the parameter~$ a $ approaches $ 0 $. In this limit, the Lie algebra $ \mathfrak{g}_{k,a} = \operatorname{span}_\mathbb{R} \bigl\{\mathbb{H}_{k,a}, \mathbb{E}^+_{k,a}, \mathbb{E}^-_{k,a}\bigr\} \cong \mathfrak{sl}(2, \mathbb{R}) $ contracts to a three-dimensional commutative Lie algebra $ \mathfrak{g}_{k,0}$, and its spectral properties change. We describe the joint spectral decomposition for $ \mathfrak{g}_{k,0}$, and discuss formulas for operator semigroups with infinitesimal generators in $ \mathfrak{g}_{k,0}$. In particular, we describe the integral kernel of $ \exp\bigl(z |x|^2 \Delta_k\bigr) $ as an infinite series, which, in some low-dimensional cases, can be expressed in a closed form using the theta function.}


\Classification{43A32; 22E45; 47D03}

\section{Introduction}

\subsection{Background}

A minimal representation is an infinite-dimensional irreducible representation
of a simple Lie group with the smallest Gelfand--Kirillov dimension. However,
at the same time, it can be thought of as a manifestation of large symmetry of
the space acted on by the group, and hence, it is expected to control global
analysis on the space effectively. This is the idea of \emph{global analysis of
minimal representations} initiated by T.~Kobayashi~\cite{MR2849643,MR3070645},
which led a transition from algebraic representation theory to analytic
representation theory. See also \cite[Section~VII]{MR4867008} for an excellent survey.

From the viewpoint of global analysis of minimal representations, the
classical Fourier transform on the Euclidean space $ \setR^N $ can be
interpreted as a unitary inversion operator in the Weil representation, which is
a unitary representation of the metaplectic group $ \mathrm{Mp}(N, \setR) $ on
the Hilbert space $ L^2\bigl(\setR^N\bigr) $ (see \cite{MR983366} for more details) and
decomposes into two irreducible components, each of which is a minimal
representation. Promoting this interpretation,
Kobayashi--Mano~\cite{MR2134314,MR2317306,MR2401813,MR2858535} introduced the \emph{Fourier
transform on the light cone} as a unitary inversion operator in an $ L^2 $-model
of a minimal representation of $ \mathrm{O}(p, q) $ and developed a new theory of
harmonic analysis. The special case $(p, q) = (N + 1, 2) $, where the model
Hilbert space is isomorphic to $ L^2\bigl(\setR^N, \enorm{x}^{-1} {\rm d}x\bigr) $, is studied
in \cite{MR2134314,MR2401813}.

After that, Ben~Sa\"{\i}d--Kobayashi--{\O}rsted~\cite{MR2566988,MR2956043} introduced a family of
$ \mathfrak{sl}_2 $-triples of differential-difference operators
$
 \DiffH{k, a}$, $ \DiffEp{k, a}$, $ \DiffEm{k, a}$
on $ \setR^N \setminus \setenum{0} $ indexed by two parameters $ k $ and $ a $,
and defined the \emph{$ (k, a) $-generalized Laguerre semigroup}
\[
 \mscrI_{k,a}(z)
 = \exp\Paren*{\frac{z}{\rm i} \bigl(\DiffEm{k, a} - \DiffEp{k, a}\bigr)},
 \qquad \RePart z \geq 0
\]
and the \emph{$ (k, a) $-generalized Fourier transform}
\smash{$\mscrF_{k,a} = {\rm e}^{\frac{{\rm i}\pi}{2} \frac{2 \dindex{k} + a + N - 2}{a}}
 \mscrI_{k,a}\Paren*{\frac{{\rm i}\pi}{2}}$}.
Here, $ k $ is a combinatorial parameter derived from the Dunkl operators,
and $ a > 0 $ is a deformation parameter. The $ (k, a) $-generalized
Fourier transform $ \mscrF_{k,a} $ includes some known transforms:
\begin{itemize}\itemsep=0pt
 \item The $ (0, 2) $-generalized Fourier transform $ \mscrF_{0, 2} $ is
 the classical Fourier transform.
 \item The $ (0, 1) $-generalized Fourier transform $ \mscrF_{0, 1} $ is
 the Hankel transform, or the Fourier transform on the light cone
 for $ (p, q) = (N + 1, 2) $.
 \item The $ (k, 2) $-generalized Fourier transform $ \mscrF_{k, 2} $ is
 the Dunkl transform~\cite{MR1199124}.
\end{itemize}
The parameter $ a $ therefore provides a continuous interpolation between
the two minimal representations of the simple Lie groups
$ \mathrm{Mp}(N, \mathbb{R}) $ and $ \mathrm{O}(N + 1, 2) $.

\subsection{Results of the paper}

Let \smash{$ \mfrakg_{k,a} = \lspan_\setR \bigl\{\DiffH{k, a}, \DiffEp{k, a}, \DiffEm{k, a}\bigr\}
\cong \mathfrak{sl}(2, \setR) $}.
Ben~Sa\"{\i}d--Kobayashi--{\O}rsted~\cite[Theorems~3.30 and 3.31]{MR2956043} showed that,
for $ a > 0 $, the action of $ \mfrakg_{k,a} $ on $ L^2(\setR^N, w_{k,a}(x) {\rm d}x) $
(see \eqref{eq:weight} for the definition of $ w_{k,a} $)
lifts to a unique unitary representation of the universal covering Lie group~\smash{$ \widetilde{\mathrm{SL}}(2, \setR) $} of $ \mathrm{SL}(2, \setR) $ and found
its irreducible decomposition explicitly; the Hilbert space $ L^2\bigl(\setR^N,\allowbreak w_{k,a}(x) {\rm d}x\bigr)$ decomposes discretely
with finite multiplicities into relatively discrete series representations
of~\smash{$ \widetilde{\mathrm{SL}}(2, \setR) $}.
Furthermore, we investigated in \cite{arXiv2407-01345} the case $ a < 0 $,
which provided an extension of the parameter~$ a $.

In the present paper, we study the behavior as $ a \to 0 $.
Although the operators $ \DiffH{k, a} $, \smash{$ \DiffEp{k, a} $} or $ \DiffEm{k, a} $
are not well-defined for $ a = 0 $, the Lie algebra
\smash{$ \mfrakg_{k,a} \cong \mathfrak{sl}(2, \setR) $} contracts
to a three-dimensional commutative Lie algebra \smash{$ \mfrakg_{k,0}\cong \setR^3 $}
as $ a \to 0 $.
Such a contraction of Lie algebras (or corresponding Lie groups) was earlier
formalized by Inonu--Wigner~\cite{MR55352}, where it is referred to as
a \emph{contraction of groups}. Classical examples include the contraction from
the orthogonal group $ \mathrm{O}(3) $ (resp.\ $ \mathrm{O}(2, 1) $) to the Euclidean motion group
$ \mathrm{O}(2) \ltimes \setR^2 $, which reflects the geometric phenomenon that the sphere of
curvature $ \kappa > 0 $ (resp.\ the hyperbolic plane of curvature $ \kappa < 0 $)
approaches the flat Euclidean plane as $ \kappa \to 0 $.

We then consider the action of $ \mfrakg_{k,0}$ on $ L^2\bigl(\setR^N, w_{k,0}(x) {\rm d}x\bigr) $
(note that the weight function $ w_{k,a} $ is well-defined even for $ a = 0 $).
As an analog of the result in the case $ a \neq 0 $, we describe that
the joint spectral decomposition for the operators in $ \mfrakg_{k,0}$ on
$ L^2\bigl(\setR^N, w_{k,0}(x) {\rm d}x\bigr) $ (Theorem~\ref{thm:spectral-decomposition})
and show that it lifts to a unique unitary representation of $ \setR^3 $
(Theorem~\ref{thm:exponential}). This is the main result of the paper.
In contrast to the case $ a \neq 0 $, this spectral decomposition involves only
the continuous spectrum.

Moreover, we discuss formulas for operator semigroups with infinitesimal
generators in $ \mfrakg_{k,0}$ (see Theorems~\ref{thm:scaling} and \ref{thm:integral-kernel}).
In particular, we describe the integral kernel of $ \exp\bigl(z \enorm{x}^2 \Laplacian_k\bigr) $
as an infinite series, which, in some low-dimensional cases, can be expressed
in a closed form using the theta function
(Theorems~\ref{thm:integral-kernel-of-1-dim}, \ref{thm:integral-kernel-of-2-dim} and \ref{thm:integral-kernel-of-4-dim}).
Although the $ (k, a) $-generalized Laguerre semigroup and
the $ (k, a) $-generalized Fourier transform are not well-defined for $ a = 0 $,
the operator semigroup $ \bigl({\rm e}^{-z} \exp\bigl(z \enorm{x}^2 \Laplacian_k\bigr)\bigr)_{\RePart z \geq 0} $
may be viewed as the ``renormalized'' $ (k, a) $-generalized Laguerre semigroup
for $ a = 0 $ (Theorem~\ref{rem:laguerre-semigroup}).
Note that explicit formulas and estimates for the integral kernels of
the $ (k, a) $-generalized Laguerre semigroup and
the $ (k, a) $-generalized Fourier transform have been extensively studied
in Ben~Sa\"{\i}d--Kobayashi--{\O}rsted~\cite[Sections~4.3--4.5 and 5.2--5.4]{MR2956043}
and subsequent papers~\cite{arXiv2407-03664,MR3759078,MR3172014,MR4832104,MR4629458,MR4871318} up to the present.
There are also unpublished results by Mano and related results by Demni~\cite{MR2859027}.

Thus, this paper analyzes representation-theoretic aspects of contraction of
Lie algebras in the framework of $ (k, a) $-generalized Fourier analysis.
We note that, recently, Benoist--Kobayashi~\cite[Theorem~1.2]{MR4653761} discovered
a relationship between \emph{limit algebras} (see Section~1.4 of their paper) of
$ \mfrakh = \Lie(H) $ in $ \mfrakg = \Lie(G) $ and $ L^2 $-analysis of $ G/H $
in the context of tempered unitary representations. It can be viewed as an
application of the notion of contraction of Lie algebras to representation theory.

\subsection{Organization of the paper}

In Section~\ref{sec:preliminaries}, we will briefly review Dunkl theory and
the differential-difference operators $ \DiffH{k, a} $, $ \DiffEp{k, a} $
and $ \DiffEm{k, a} $ introduced by Ben~Sa\"{\i}d--Kobayashi--{\O}rsted.
In Section~\ref{sec:contraction}, we will discuss the contraction of the
$ \mathfrak{sl}_2 $-triple as $ a \to 0 $.
In Section~\ref{sec:closed-form}, we will give a closed-form expression for
the integral kernel of $ \exp\bigl(z \enorm{x}^2 \Laplacian_k\bigr) $
in some low-dimensional cases.

\subsection{Notation}

\begin{itemize}\itemsep=0pt
 \item $ \setN = \setenum{0, 1, 2, \dots} $.
 \item We write $ \innprod{\blank}{\blank} $ for the Euclidean inner product,
 and $ \enorm{\blank} $ for the Euclidean norm.
 \item $ \Sphere{N - 1} = \bigl\{x \in \setR^N\mid \enorm{x} = 1\bigr\} $.
 \item Function spaces, such as $ C^\infty $ spaces and $ L^2 $ spaces,
 are understood to consist of complex-valued functions.
 \item We write \smash{$ E_x = \sum_{j = 1}^{N} x_j \Pdif{x_j} $} for the Euler operator
 on $ \setR^N $, and $ E_r = r \frac{\rm d}{{\rm d}r} $ for the Euler operator on $ \setRp $.
\end{itemize}

\section{Preliminaries}
\label{sec:preliminaries}

In this section, we review Dunkl theory and the differential-difference operators
$ \DiffH{k, a} $, $ \DiffEp{k, a} $ and~\smash{$ \DiffEm{k, a} $} introduced
by Ben~Sa\"{\i}d--Kobayashi--{\O}rsted to the extent necessary for later use.
This section contains no new results.

\subsection{The Dunkl Laplacian}

Throughout this paper, we fix a reduced root system $ \mscrR $ on $ \setR^N $.
That is, we suppose that $ \mscrR $ satisfies the following conditions:
\begin{itemize}\itemsep=0pt
 \item $ \mscrR $ is a finite subset of $ \setR^N \setminus \setenum{0} $,
 \item $ \mscrR $ is stable under the orthogonal reflection $ r_\alpha $
 with respect to the hyperplane $ (\setR \alpha)^\perp $ for all
 $ \alpha \in \mscrR $, and
 \item $ \mscrR \cap \setR \alpha = \setenum{\alpha, -\alpha} $ for all
 $ \alpha \in \mscrR $.
\end{itemize}
Note that we do not impose crystallographic conditions on roots and do not
require that $ \mscrR $ spans $ \setR^N $.

The subgroup of $ \mathrm{O}(N) $ generated by all the reflections $ r_\alpha $
is called the \emph{reflection group associated with $ \mscrR $}.
We say that a function $ \map{k}{\mscrR}{\setC} $ is a \emph{multiplicity
function} if it is invariant under the natural action of the reflection group.
We usually write $ k_\alpha $ instead of $ k(\alpha) $. We say that
a multiplicity function $ k $ is \emph{non-negative} if $ k_\alpha \geq 0 $
for all $ \alpha \in \mscrR $. The \emph{index} of a multiplicity function
$ k $ is defined as
\[
 \dindex{k}
 = \frac{1}{2} \sum_{\alpha \in \mscrR} k_\alpha
 = \sum_{\alpha \in \mscrR^+} k_\alpha,
\]
where $ \mscrR^+ $ is any positive system of $ \mscrR $.

For a (not necessarily non-negative) multiplicity function $ k $,
the \emph{Dunkl Laplacian} $ \Laplacian_k $
(see \cite{MR917849} and \cite[Definition~1.1]{MR951883}) is defined by
\[
 \Laplacian_k F(x)
 = \Laplacian F(x)
 + \sum_{\alpha \in \mscrR^+} k_\alpha
 \Paren*{
 \frac{2 \innprod{\alpha}{\nabla F(x)}}{\innprod{\alpha}{x}}
 - \enorm{\alpha}^2 \frac{F(x) - F(r_\alpha(x))}{\innprod{\alpha}{x}^2}
 },
\]
where \smash{$ \Laplacian = \sum_{j = 1}^{N} \bigl(\Pdif{x_j}\bigr)^2 $} is the classical Laplacian
and $ \nabla = \bigl(\Pdif{x_1}, \dots, \Pdif{x_N}\bigr) $ is the classical gradient operator.
When $ k = 0 $, the Dunkl Laplacian $ \Laplacian_k $ reduces to the classical
Laplacian $ \Laplacian $.

Let $ \mcalP\bigl(\setR^N\bigr) $ denote the space of polynomials on $ \setR^N $ and
$ \mcalP^m\bigl(\setR^N\bigr) $ denote its subspace of homogeneous polynomials of degree $ m $.
The space of \emph{$ k $-harmonic polynomials of degree $ m $}
(see~\cite[Definition~1.5]{MR917849}) is defined as
\[
 \mcalH_k^m\bigl(\setR^N\bigr)
 = \bigl\{p \in \mcalP^m\bigl(\setR^N\bigr)\mid \Laplacian_k p = 0\bigr\},
\]
and the space of \emph{$ k $-spherical harmonics of degree $ m $} is defined as
\[
 \mcalH_k^m\bigl(\Sphere{N - 1}\bigr)
 = \bigl\{\restr{p}{\Sphere{N - 1}}\mid p \in \mcalH_k^m\bigl(\setR^N\bigr)\bigr\}.
\]
When $ k = 0 $, these are reduced to the space $ \mcalH^m\bigl(\setR^N\bigr) $
of classical harmonic polynomials and the space $ \mcalH^m\bigl(\Sphere{N - 1}\bigr) $
of classical spherical harmonics, respectively.

We write $ \mcalH\bigl(\Sphere{N - 1}\bigr) = \bigl\{\restr{p}{\Sphere{N - 1}}\mid p \in \mcalP\bigl(\setR^N\bigr)\bigr\} $.
The following fact is a generalization of the decomposition of
$ \mcalH\bigl(\Sphere{N - 1}\bigr) $ and $ L^2\bigl(\Sphere{N - 1}\bigr) $ into the spaces
$ \mcalH^m\bigl(\Sphere{N - 1}\bigr) $ of classical spherical harmonics.

\begin{Fact}[{\cite[pp.~37--39]{MR917849}}]\label{thm:decomposition}
 For a non-negative multiplicity function~$k$,
 we have the direct sum decomposition
 \[
 \mcalH\bigl(\Sphere{N - 1}\bigr)
 = \bigoplus_{m \in \setN} \mcalH_k^m\bigl(\Sphere{N - 1}\bigr)
 \]
 and the orthogonal decomposition
 \[
 L^2\bigl(\Sphere{N - 1}, w_k(\omega) {\rm d}\omega\bigr)
 = \sumoplus_{m \in \setN} \mcalH_k^m\bigl(\Sphere{N - 1}\bigr),
 \]
 where the weight function $ w_k $ with respect to the standard measure
 $ {\rm d}\omega $ on $ \Sphere{N - 1} $ is defined by
 \[
 w_k(\omega)
 = \prod_{\alpha \in \mscrR} \abs{\innprod{\alpha}{\omega}}^{k_\alpha}
 = \prod_{\alpha \in \mscrR^+} \abs{\innprod{\alpha}{\omega}}^{2k_\alpha}.
 \]
\end{Fact}

\subsection{Dunkl's intertwining operator and the Poisson kernel}

For a non-negative multiplicity function $ k $, Dunkl introduced a linear
automorphism $ V_k $ of the space $ \mcalP\bigl(\setR^N\bigr) $ of polynomials
that satisfies a certain intertwining property (\emph{Dunkl's intertwining operator}).
See \cite[Definition~2.2 and Theorem~2.3]{MR1145585} for the definition
and a characterization of~$ V_k $. We note that \smash{$ V_0 = \id_{\mcalP(\setR^N)} $}.

The following integral representation of Dunkl's intertwining operator $ V_k $
was obtained by R\"osler.

\begin{Fact}[{\cite[Theorem~1.2]{MR1695797}}]\label{thm:dunkls-intertwining-operator-as-integral}
 Let $ k $ be a non-negative multiplicity function.
 For each $ x \in \setR^N $, there exists a unique probability Borel measure
 $ \mu_{k, x} $ on $ \setR^N $ such that
 \[
 V_k p(x)
 = \int_{\setR^N} p(\xi) \, {\rm d}\mu_{k, x}(\xi)
 \]
 for all $ p \in \mcalP\bigl(\setR^N\bigr) $. Moreover, the support of $ \mu_{k, x} $ is
 contained in the ball $ \bigl\{\xi \in \setR^N\mid \enorm{\xi} \leq \enorm{x}\bigr\} $,
 and we have
$
 \mu_{k, x}(S)
 = \mu_{k, gx}(gS)
 = \mu_{k, rx}(rS)
$
 for any element $ g $ of the reflection group, $ r > 0 $,
 and Borel set $ S \subseteq \setR^N $.
\end{Fact}

Let $ k $ be a non-negative multiplicity function and $ m \in \setN $.
We consider the orthogonal projection \smash{$ \Pi_k^{(m)} $}
from $ L^2\bigl(\Sphere{N - 1}, w_k(\omega) {\rm d}\omega\bigr) $
onto $ \mcalH_k^m\bigl(\Sphere{N - 1}\bigr) $ and its normalized integral kernel~\smash{$ P_k^{(m)} $},
which is called the \emph{Poisson kernel} of the space of $ k $-spherical harmonics
of degree $ m $. That is, the function \smash{$ P_k^{(m)} $}
on $ \Sphere{N - 1} \times \Sphere{N - 1} $ is characterized by the formula
\begin{equation}
 \Pi_k^{(m)} p(\omega)
 = \frac{1}{\vol_k\bigl(\Sphere{N - 1}\bigr)}
 \int_{\Sphere{N - 1}} P_k^{(m)}\big(\omega, \omega'\big) p\big(\omega'\big) w_k\big(\omega'\big) \, {\rm d}\omega'
 \label{eq:poisson-kernel}
\end{equation}
for all $ p \in L^2\bigl(\Sphere{N - 1}, w_k(\omega) {\rm d}\omega\bigr) $, where
\begin{gather}\label{eq:volume}
 \vol_k\bigl(\Sphere{N - 1}\bigr)
 = \int_{\Sphere{N - 1}} \, {\rm d}w_k(\omega).
\end{gather}
Equivalently, the Poisson kernel \smash{$ P_k^{(m)} $} is given by
\[
 P_k^{(m)}(\omega, \omega')
 = \vol_k\bigl(\Sphere{N - 1}\bigr) \sum_{j = 1}^{d} p_j(\omega) p_j(\omega'),
\]
where $ (p_1, \dots, p_d) $ is an orthonormal basis of $ \mcalH_k^m\bigl(\Sphere{N - 1}\bigr) $,
regarded as a subspace of $ L^2\bigl(\Sphere{N - 1},\allowbreak w_k(\omega) {\rm d}\omega\bigr) $.

The Poisson kernel \smash{$ P_k^{(m)} $} can be expressed in terms of Dunkl's intertwining
operator and the Gegenbauer polynomials. To state this result, we first prepare
some notation. For $ \nu \in \setC $ and~${ m \in \setN }$,
we consider the Gegenbauer polynomial $ C_m^\nu $ defined by the generating
function
\begin{equation}
 \bigl(1 - 2t\xi + \xi^2\bigr)^{-\nu}
 = \sum_{m = 0}^{\infty} C_m^\nu(t) \xi^m,
 \label{eq:gegenbauer}
\end{equation}
and the renormalized Gegenbauer polynomial $ \widetilde{C}_m^\nu $ defined by
\begin{equation}
 \widetilde{C}_m^\nu(t)
 = \frac{m + \nu}{\nu} C_m^\nu(t).
 \label{eq:renormalized-gegenbauer}
\end{equation}
For $ \nu = 0 $, we define $ \widetilde{C}_m^0 $ by the limit formula
(see \cite[equation~(6.4.13)]{MR1688958})
\begin{equation}
 \widetilde{C}_m^0(t)
 = \lim_{\nu \to 0} \widetilde{C}_m^\nu(t)
 = \begin{cases}
 1, & m = 0 , \\
 2 T_m(t), & m \geq 1,
 \end{cases}
 \label{eq:limit-formula}
\end{equation}
where $ T_m $ denotes the Chebyshev polynomial of the first kind, which is
characterized by the formula $ T_m(\cos \theta) = \cos m\theta $.

\begin{Fact}[{\cite[Theorem~4.1]{MR1145585}}]\label{thm:poisson-kernel}
 Let $ k $ be a non-negative multiplicity function and $ m \in \setN $.
 The Poisson kernel \smash{$ P_k^{(m)} $} is given by
 \[
 P_k^{(m)}(\omega, \omega')
 = V_k\big(\widetilde{C}_m^{\frac{2 \dindex{k} + N - 2}{2}}\big(\innprod{\blank}{\omega'}\big)\big)(\omega).
 \]
\end{Fact}

\begin{Remark}\label{rem:classical-poisson-kernel}
 When $ k = 0 $, we have $ V_0 = \id_{\mcalP(\setR^N)} $, so that
 the formula in Theorem~\ref{thm:poisson-kernel} reduces~to
 \[
 P_0^{(m)}\bigl(\omega, \omega'\bigr)
 = \widetilde{C}_m^{\frac{N - 2}{2}}\bigl(\innprod{\omega}{\omega'}\bigr).
 \]
 See \cite[Theorem~9.6.3 and Remark~9.6.1]{MR1688958}.
\end{Remark}

\begin{Remark}\label{rem:poisson-kernel-of-1-dim}
 In the case $ N = 1 $, we have $ \Sphere{0} = \setenum{\pm 1} $ and
 the formula in Theorem~\ref{thm:poisson-kernel} reduces~to
 \[
 P_k^{(m)}(\omega, \omega')
 = \begin{cases}
 1, & m = 0 ,\\
 \sgn(\omega\omega') ,& m = 1, \\
 0, & m \geq 2,
 \end{cases}
 \]
 which corresponds to the fact that
 \[
 \mcalH_k^m\big(\Sphere{0}\big)
 = \begin{cases}
 \setC 1, & m = 0, \\
 \setC \sgn, & m = 1 ,\\
 0, & m \geq 2.
 \end{cases}
 \]
 Here, $ 1 $ and $ \sgn $ denote the constant function and the sign function on
 $ \Sphere{0} = \setenum{\pm 1} $, respectively.
\end{Remark}

\subsection[The differential-difference operators H\_\{k,a\}, E\^{}+\_\{k,a\} and E\^{}-\_\{k,a\}]{The differential-difference operators $\boldsymbol{\DiffH{k,a}}$, $ \boldsymbol{\DiffEp{k,a}}$ and $ \boldsymbol{\DiffEm{k,a}}$}

Let $ k $ be a multiplicity function and $ a \in \setC \setminus \setenum{0} $.
We recall the definition of the differential-difference operators
$ \DiffH{k, a} $, $ \DiffEp{k, a} $ and $ \DiffEm{k, a} $
on $ \setR^N \setminus \setenum{0} $ from \cite[equation~(3.3)]{MR2956043}:
\[
 \DiffH{k, a}
 = \frac{2}{a} E_x + \frac{2 \dindex{k} + a + N - 2}{a}, \qquad
 \DiffEp{k, a}
 = \frac{\rm i}{a} \enorm{x}^a, \qquad
 \DiffEm{k, a}
 = \frac{\rm i}{a} \enorm{x}^{2 - a} \Laplacian_k.
\]

Additionally, for $ m \in \setN $, we consider the following differential
operators on $ \setRp $:
\begin{gather*}
 \DiffH{k, a}[m]
 = \frac{2}{a} E_r + \frac{2 \dindex{k} + a + N - 2}{a}, \qquad
 \DiffEp{k, a}[m]
 = \frac{\rm i}{a} r^a, \\
 \DiffEm{k, a}[m]
 = \frac{\rm i}{a} r^{-a} (E_r - m)(E_r + m + 2 \dindex{k} + N - 2).
\end{gather*}
These are the radial parts of $ \DiffH{k, a} $, $ \DiffEp{k, a} $ and
$ \DiffEm{k, a} $ respectively in the following sense.

\begin{Proposition}\label{thm:radial-parts}
 Let $ k $ be a multiplicity function, $ a \in \setC \setminus \setenum{0} $,
 and $ m \in \setN $.
 For $ p \in \mcalH_k^m\bigl(\Sphere{N - 1}\bigr) $ and $ f \in C^\infty(\setRp) $, we have
 \begin{gather*}
 \DiffH{k, a} (p \otimes f) = p \otimes \DiffH{k, a}[m] f, \qquad
 \DiffEp{k, a} (p \otimes f) = p \otimes \DiffEp{k, a}[m] f, \\
 \DiffEm{k, a} (p \otimes f) = p \otimes \DiffEm{k, a}[m] f,
 \end{gather*}
 where $ p \otimes f $ denotes the function $ r\omega \mapsto p(\omega) f(r) $
 on $ \setR^N \setminus \setenum{0} $.
\end{Proposition}

\begin{proof}
 The first and second equations are clear. We now prove the third equation.
 We use the polar coordinates $ x = r\omega $, where $ r \in \setRp $ and
 $ \omega \in \Sphere{N - 1} $. We extend $ p \in \mcalH_k^m\bigl(\Sphere{N - 1}\bigr) $
 to a~$ k $-harmonic polynomial of degree $ m $ on $ \setR^N $, which we again
 write $ p $. Then, since $ \Laplacian_k p = 0 $, we have\looseness=-1
 \begin{equation}
 \Laplacian_k (p \otimes f)(x)
 = \Laplacian_k (r^{-m} f(r) p(x))
 = [\Laplacian_k, r^{-m} f(r)] p(x),
 \label{eq:radial-parts-1}
 \end{equation}
 where $ [\blank, \blank] $ denotes the commutator of operators.
 For a radial function $ g(r) $, the commutator~${ [\Laplacian_k, g(r)] }$
 can be computed by the Leibniz rule as
 \begin{align*}
 [\Laplacian_k, g(r)]
 &= \Laplacian (g(r))
 + 2 \innprod{\nabla (g(r))}{\nabla}
 + \sum_{\alpha \in \mscrR^+} k_\alpha \frac{2 \innprod{\alpha}{\nabla (g(r))}}{\innprod{\alpha}{x}} \\
 &= g''(r) + \frac{1}{r} g'(r) (2E_x + 2 \dindex{k} + N - 1) \\
 &= r^{-2} \bigl(E_r^2 g(r) + E_r g(r) (2E_x + 2 \dindex{k} + N - 2)\bigr).
 \end{align*}
 Setting $ g(r) = r^{-m} f(r) $ and applying this commutator to $ p(x) $,
 we have
 \begin{gather}
 [\Laplacian_k, r^{-m} f(r)] p(x)
 \nonumber \\
 \qquad= r^{-2} \bigl(E_r^2 (r^{-m} f(r)) + E_r (r^{-m} f(r)) (2E_x + 2 \dindex{k} + N - 2)\bigr) p(x)
 \nonumber \\
 \qquad= r^{-m - 2} \bigl((E_r - m)^2 f(r) + (E_r - m) f(r) (2E_x + 2 \dindex{k} + N - 2)\bigr) p(x)
 \nonumber \\
 \qquad= r^{-m - 2} \bigl((E_r - m)^2 + (E_r - m) (2m + 2 \dindex{k} + N - 2)\bigr) f(r) p(x)
 \nonumber \\
 \qquad= r^{-2} (E_r - m) (E_r + m + 2 \dindex{k} + N - 2) f(r) p(\omega)
 \nonumber \\
 \qquad= \bigl(p \otimes r^{-2} (E_r - m) (E_r + m + 2 \dindex{k} + N - 2) f\bigr)(x).
 \label{eq:radial-parts-2}
 \end{gather}
 The third equation follows from \eqref{eq:radial-parts-1} and \eqref{eq:radial-parts-2}.
\end{proof}

\begin{Proposition}\label{thm:sl2-triple}
 Let $ k $ be a multiplicity function and $ a \in \setC \setminus \setenum{0} $.
 \begin{enumarabicp}
 \item The differential-difference operators
 $ \DiffH{k, a} $, $ \DiffEp{k, a} $ and $ \DiffEm{k, a} $
 form an $ \mathfrak{sl}_2 $-triple. That is,
 \[
 \big[\DiffH{k, a}, \DiffEp{k, a}\big] = 2 \DiffEp{k, a}, \qquad
 [\DiffH{k, a}, \DiffEm{k, a}] = -2 \DiffEm{k, a}, \qquad
 \big[\DiffEp{k, a}, \DiffEm{k, a}\big] = \DiffH{k, a}.
 \]
 \item For any $ m \in \setN $, the differential operators
 $ \DiffH{k, a}[m] $, $ \DiffEp{k, a}[m] $ and $ \DiffEm{k, a}[m] $
 form an $ \mathfrak{sl}_2 $-triple. That~is,
 \[
 \big[\DiffH{k, a}[m], \DiffEp{k, a}[m]\big] = 2 \DiffEp{k, a}, \qquad
 \big[\DiffH{k, a}[m], \DiffEm{k, a}[m]\big] = -2 \DiffEm{k, a}, \qquad
 \big[\DiffEp{k, a}[m], \DiffEm{k, a}[m]\big] = \DiffH{k, a}.
 \]
 \end{enumarabicp}
\end{Proposition}

\begin{proof}
(1)
 It is \cite[Theorem~3.2]{MR2956043}.
(2)
 It follows from (1) and Theorem~\ref{thm:radial-parts}.
\end{proof}

\section[Contraction of the sl\_2-triple as a to 0]{Contraction of the $\boldsymbol{\mathfrak{sl}_2}$-triple as $\boldsymbol{a \to 0}$}
\label{sec:contraction}

\subsection[The commutative Lie algebras g\_\{k,0\} and g\_\{k,0\}\^{}\{rad\}]{The commutative Lie algebras $\boldsymbol{\mfrakg_{k,0}}$ and $ \boldsymbol{\mfrakg_{k,0}^{\mathrm{rad}}}$}

For a multiplicity function $ k $ and $ a \in \setC \setminus \setenum{0} $,
we write
\begin{align*}
 \mfrakg_{k,a}
 &= \lspan_\setR \bigl\{\DiffH{k, a}, \DiffEp{k, a}, \DiffEm{k, a}\bigr\} = \lspan_\setR \bigl\{a \DiffH{k, a}, a \DiffEp{k, a}, a \DiffEm{k, a}\bigr\} \\
 &= \lspan_\setR \bigl\{2E_x + 2 \dindex{k} + a + N - 2, {\rm i} \enorm{x}^a,{\rm i} \enorm{x}^{2 - a} \Laplacian_k\bigr\}.
\end{align*}
Putting $ a = 0 $ in the above equation, we define
\begin{equation}
 \mfrakg_{k,0}
 = \lspan_\setR \bigl\{2E_x + 2 \dindex{k} + N - 2, {\rm i}, {\rm i} \enorm{x}^2 \Laplacian_k\bigr\}.
 \label{eq:gk0}
\end{equation}
Similarly, we write
\begin{align*}
 \mfrakg_{k,a}^{(m)}
 &= \lspan_\setR \bigl\{\DiffH{k, a}[m], \DiffEp{k, a}[m], \DiffEm{k, a}[m]\bigr\}= \lspan_\setR \bigl\{a \DiffH{k, a}[m], a \DiffEp{k, a}[m], a \DiffEm{k, a}[m]\bigr\} \\
 &= \lspan_\setR \bigl\{2E_r + 2 \dindex{k} + a + N - 2, {\rm i}r^a, {\rm i}r^{-a} (E_r - m)(E_r + m + 2 \dindex{k} + N - 2)\bigr\}.
\end{align*}
and define
\begin{align}
 \mfrakg_{k,0}^{\mathrm{rad}}
 &= \lspan_\setR \{2E_r + 2 \dindex{k} + N - 2, {\rm i}, {\rm i} (E_r - m)(E_r + m + 2 \dindex{k} + N - 2)\} \nonumber \\
 &= \lspan_\setR \bigl\{2E_r + 2 \dindex{k} + N - 2, {\rm i}, {\rm i} \bigl(E_r^2 + (2 \dindex{k} + N - 2) E_r - m (m + 2 \dindex{k} + N - 2)\bigr)\bigr\} \nonumber \\
 &= \lspan_\setR \bigl\{2E_r + 2 \dindex{k} + N - 2, {\rm i}, {\rm i} \bigl(E_r^2 + (2 \dindex{k} + N - 2) E_r\bigr)\bigr\}.
 \label{eq:gk0rad}
\end{align}
Note that the right-hand side of the above definition does not depend on $ m $,
which justifies the notation $ \mfrakg_{k,0}^{\mathrm{rad}} $.

\begin{Proposition}\label{thm:radial-parts-zero}
 Let $ k $ be a multiplicity function.
 For $ p \in \mcalH_k^m\bigl(\Sphere{N - 1}\bigr) $ and $ f \in C^\infty(\setRp) $, we have
 \begin{gather*}
 (2E_x + 2 \dindex{k} + N - 2) (p \otimes f)
 = p \otimes (2E_r + 2 \dindex{k} + N - 2) f, \\
{\rm i} (p \otimes f)
 = p \otimes {\rm i}f, \qquad
{\rm i} \enorm{x}^2 \Laplacian_k (p \otimes f)
 = p \otimes {\rm i} (E_r - m)(E_r + m + 2 \dindex{k} + N - 2) f,
 \end{gather*}
 where $ p \otimes f $ denotes the function $ r\omega \mapsto p(\omega) f(r) $
 on $ \setR^N \setminus \setenum{0} $.
\end{Proposition}

\begin{proof}
 The proof goes along the same lines as that of Theorem~\ref{thm:radial-parts}.
\end{proof}

\begin{Proposition}
 Let $ k $ be a multiplicity function.
 \begin{enumarabicp}
 \item The space $ \mfrakg_{k,0}$ of differential-difference operators on
 $ \setR^N \setminus \setenum{0} $ is a three-dimensional commutative Lie algebra.
 \item The space $ \mfrakg_{k,0}^{\mathrm{rad}} $ of differential operators on
 $ \setRp $ is a three-dimensional commutative Lie algebra.
 \end{enumarabicp}
\end{Proposition}

\begin{proof}
(1)
 By Theorem~\ref{thm:sl2-triple}, we have
 \[
 \big[a \DiffH{k, a}, a \DiffEp{k, a}\big] = 2a \cdot a \DiffEp{k, a}, \qquad
 [a \DiffH{k, a}, a \DiffEm{k, a}] = -2a \cdot a \DiffEm{k, a}, \qquad
 \big[a \DiffEp{k, a}, a \DiffEm{k, a}\big] = a \cdot a \DiffH{k, a}.
 \]
 By taking the limit as $ a \to 0 $, we have
 \[
 [2E_x + 2 \dindex{k} + N - 2, {\rm i}] = 0, \qquad
 \big[2E_x + 2 \dindex{k} + N - 2, {\rm i} \enorm{x}^2 \Laplacian_k\big] = 0, \qquad
 \big[{\rm i}, {\rm i} \enorm{x}^2 \Laplacian_k\big] = 0.
 \]
 (It can also be shown by a direct computation.)

(2) It follows from (1) and Theorem~\ref{thm:radial-parts-zero}.
\end{proof}

\subsection[Joint spectral decomposition for g\_\{k,0\} and g\_\{k,0\}\^\{rad\}]{Joint spectral decomposition for $ \boldsymbol{\mfrakg_{k,0}} $ and $\boldsymbol{\mfrakg_{k,0}^{\mathrm{rad}}}$}

In the following, we consider a \emph{non-negative} multiplicity function $ k $.
In the next two theorems, we use the unitary operator
\begin{gather*}
 U_{N, k}\colon\ L^2\big(\setRp, r^{2 \dindex{k} + N - 3} {\rm d}r\big)\to {L^2(\setR, {\rm d}s)}, \\
 U_{N, k} f(s) = {\rm e}^{\frac{2 \dindex{k} + N - 2}{2} s} f({\rm e}^s), \qquad
 U_{N, k}^{-1} g(r) = r^{-\frac{2 \dindex{k} + N - 2}{2}} g(\log r)
\end{gather*}
and the (classical) Fourier transform
\begin{gather*}
 \mscrF\colon \ L^2(\setR, {\rm d}s)\to {L^2(\setR, {\rm d}\sigma)}, \\
 \mscrF g(\sigma) = \frac{1}{\sqrt{2\pi}} \int_\setR g(s) {\rm e}^{-{\rm i}\sigma s} \, {\rm d}s, \qquad
 \mscrF^{-1} h(s) = \frac{1}{\sqrt{2\pi}} \int_\setR g(s) {\rm e}^{{\rm i}\sigma s} \, {\rm d}s.
\end{gather*}

We recall some terminology related to operators on a Hilbert space.
A densely defined operator $ T $ on a Hilbert space is called
\emph{self-adjoint} (resp.\ \emph{skew-adjoint}) if its adjoint $ T^* $
is equal to $ T $ (resp.\ ${\rm i}T $), that is, these have the same domain and
coincide on it. A closable operator~$ T $ on a Hilbert space is called
\emph{essentially self-adjoint} (resp.\ \emph{essentially skew-adjoint})
if its closure~$ \closure{T} $ is self-adjoint (resp.\ skew-adjoint).

\begin{Theorem}\label{thm:spectral-decomposition-rad}
 Let $ k $ be a non-negative multiplicity function.
 Every differential operator in~$ \mfrakg_{k,0}^{\mathrm{rad}} $
 $($see {\rm\eqref{eq:gk0rad}} for the definition$)$ defined on
 the domain $ C_{\mathrm{c}}^\infty(\setRp) $ is an essentially skew-adjoint
 operator on $
 L^2\bigl(\setRp, r^{2 \dindex{k} + N - 3} {\rm d}r\bigr) $.
 Moreover, via the unitary operator
 \[ \mscrF \circ U_{N, k}\colon \ L^2\bigl(\setRp, r^{2 \dindex{k} + N - 3} {\rm d}r\bigr)\to L^2(\setR, {\rm d}\sigma) ,\]
 the closures of
 \begin{gather*}
 \restr{(2E_r + 2 \dindex{k} + N - 2)}{C_{\mathrm{c}}^\infty(\setRp)}, \qquad
 {\rm i}\, \id_{C_{\mathrm{c}}^\infty(\setRp)}, \\
 \restr{{\rm i} (E_r - m) (E_r + m + 2 \dindex{k} + N - 2)}{C_{\mathrm{c}}^\infty(\setRp)}
 \end{gather*}
 correspond to the multiplication operators
 \[
 2{\rm i}\sigma, \qquad
{\rm i}, \qquad
 -{\rm i} \Paren*{\sigma^2 + \Paren*{m + \frac{2 \dindex{k} + N - 2}{2}}^2},
 \]
 respectively.
\end{Theorem}

\begin{proof}
 Via the unitary operator $ \map{U_{N, k}}{L^2\bigl(\setRp, r^{2 \dindex{k} + N - 3} {\rm d}r\bigr)}{L^2(\setR, {\rm d}s)} $,
 the operator \[
 \left(E_r + \frac{2 \dindex{k} + N - 2}{2}\right)\bigg|_{C_{\mathrm{c}}^\infty(\setR)} \]
 corresponds to \smash{$ \restr{\frac{\rm d}{{\rm d}s}}{C_{\mathrm{c}}^\infty(\setR)} $}.
 As is well-known, for any complex polynomial $ P(\sigma) $ such that
 $ P({\rm i}\sigma) $ is real-valued, $ \restr{P\bigl(\frac{\rm d}{{\rm d}s}\bigr)}{C_{\mathrm{c}}^\infty(\setR)} $
 is an essentially self-adjoint operator on $ L^2(\setR,{\rm d}s) $ and
 its closure corresponds to the multiplication operator by the function
 $ P({\rm i}\sigma) $ via the Fourier transform $ \map{\mscrF}{L^2(\setR, {\rm d}s)}{L^2(\setR, {\rm d}\sigma)} $.
 Since $ 2E_r + 2 \dindex{k} + N - 2 $, ${\rm i} $ and $ {\rm i} (E_r - m) (E_r + m + 2 \dindex{k} + N - 2) $
 can be expressed as $ {\rm i} $ times such polynomials of \smash{$ E_r + \frac{2 \dindex{k} + N - 2}{2} $},
 now the assertion follows.
\end{proof}

We recall that the $ L^2 $-theory of the $ \mathfrak{sl}_2 $-triple
\smash{$ \bigl(\DiffH{k, a}, \DiffEp{k, a}, \DiffEm{k, a}\bigr) $} was considered on
the Hilbert space \smash{$ L^2\bigl(\setR^N, w_{k,a}(x) {\rm d}x\bigr) $}, where the weight function
$ \map{w_{k,a}}{\setR^N}{\setRzp} $
(see \cite[equation~(1.2)]{MR2956043}, $ \vartheta_{k,a} $~in their notation) is defined by
\begin{equation}
 w_{k,a}(x)
 = \enorm{x}^{a - 2} \prod_{\alpha \in \mscrR} \abs{\innprod{\alpha}{x}}^{k_\alpha}
 = \enorm{x}^{a - 2} \prod_{\alpha \in \mscrR^+} \abs{\innprod{\alpha}{x}}^{2k_\alpha}.
 \label{eq:weight}
\end{equation}
By the polar decomposition
\smash{$ w_{k,a}(x) {\rm d}x = w_k(\omega) {\rm d}\omega \otimes r^{2 \dindex{k} + a + N - 3} {\rm d}r $}
(here, $ w_k $ is as defined in Theorem~\ref{thm:decomposition}) and
Theorem~\ref{thm:decomposition}, we have the orthogonal decomposition
\begin{align}
 L^2\bigl(\setR^N, w_{k,a}(x) {\rm d}x\bigr)
 &= L^2\bigl(\Sphere{N - 1}, w_k(\omega) {\rm d}\omega\bigr) \otimeshat L^2\bigl(\setRp, r^{2 \dindex{k} + a + N - 3} {\rm d}r\bigr)
 \nonumber \\
 &= \sumoplus_{m \in \setN} \mcalH_k^m\bigl(\Sphere{N - 1}\bigr) \otimes L^2\bigl(\setRp, r^{2 \dindex{k} + a + N - 3} {\rm d}r\bigr).
 \label{eq:decomposition}
\end{align}
Note that $ w_{k,a} $ is well-defined and the above orthogonal decomposition
holds even for $ a = 0 $.

We now state the main result of this paper.

\begin{Theorem}\label{thm:spectral-decomposition}
 Let $ k $ be a non-negative multiplicity function.
 Every differential-difference operator in $ \mfrakg_{k,0}$
 $($see~\eqref{eq:gk0} for the definition$)$ defined on the domain
 \[
 \mcalD
 = \mcalH\bigl(\Sphere{N - 1}\bigr) \otimes C_{\mathrm{c}}^\infty(\setRp)
 = \lspan_\setC \bigl\{p \otimes f\mid p \in \mcalH\bigl(\Sphere{N - 1}\bigr), \, f \in C_{\mathrm{c}}^\infty(\setRp) \bigr\},
 \]
 is an essentially skew-adjoint operator on $ L^2\bigl(\setR^N, w_{k,0}(x) {\rm d}x\bigr) $.
 Moreover, via the unitary operator
 \[
\id_{L^2(\Sphere{N - 1}, w_k(\omega) {\rm d}\omega)} \otimeshat (\mscrF \circ U_{N, k})\colon\
 L^2\bigl(\setR^N, w_{k,0}(x) {\rm d}x\bigr)\to
 L^2\bigl(\Sphere{N - 1}, w_k(\omega) {\rm d}\omega\bigr) \otimeshat L^2(\setR, {\rm d}\sigma),
 \]
 the closures of
 \[
 \restr{(2E_x + 2 \dindex{k} + N - 2)}{\mcalD}, \qquad
 {\rm i}\, \id_\mcalD, \qquad
 \restr{\big( {\rm i}\enorm{x}^2 \Laplacian_k\big)}{\mcalD}
 \]
 correspond to the multiplication operators
 \[
 \id_{L^2(\Sphere{N - 1}, w_k(\omega) {\rm d}\omega)} \otimeshat 2{\rm i}\sigma, \qquad
 {\rm i}, \qquad
 \sumoplus_{m \in \setN} \id_{\mcalH_k^m(\Sphere{N - 1})} \otimes \Paren*{- {\rm i} \Paren*{\sigma^2 + \Paren*{m + \frac{2 \dindex{k} + N - 2}{2}}^2}},
 \]
 respectively.
\end{Theorem}

\begin{proof}
 It follows from Theorems~\ref{thm:radial-parts-zero} and~\ref{thm:spectral-decomposition-rad},
 and the orthogonal decomposition \eqref{eq:decomposition}.
\end{proof}

\begin{Remark}
 Since the unitary operator
 $ \map{\mscrF \circ U_{N, k}}{L^2(\setRp, r^{2 \dindex{k} + N - 3} {\rm d}r)}{L^2(\setR, {\rm d}\sigma)} $
 ``maps'' \smash{$ \frac{1}{\sqrt{2\pi}} r^{-\frac{2 \dindex{k} + N - 2}{2} + {\rm i}\sigma} $} to the Dirac
 distribution $ \delta_\sigma $, we have the direct integral decomposition
 \[
 L^2\bigl(\setRp, r^{2 \dindex{k} + N - 3} {\rm d}r\bigr)
 = \intoplus_\setR \setC \frac{1}{\sqrt{2\pi}} r^{-\frac{2 \dindex{k} + N - 2}{2} + {\rm i}\sigma} \, {\rm d}\sigma
 \]
 and may write Theorem~\ref{thm:spectral-decomposition-rad} as
 \begin{gather*}
 2E_r + 2 \dindex{k} + N - 2
 = \intoplus_\setR 2{\rm i}\sigma\, {\rm d}\sigma, \qquad
 {\rm i}\, \id
 = \intoplus_\setR {\rm i} \, {\rm d}\sigma, \\
{\rm i} (E_r - m) (E_r + m + 2 \dindex{k} + N - 2)
 = \intoplus_\setR \Paren*{-{\rm i} \Paren*{\sigma^2 + \Paren*{m + \frac{2 \dindex{k} + N - 2}{2}}^2}}\, {\rm d}\sigma.
 \end{gather*}
 Similarly, we have the direct integral decomposition
 \[
 L^2\bigl(\setR^N, w_{k,0}(x) {\rm d}x\bigr)
 = \sumoplus_{m \in \setN} \intoplus_\setR
 \mcalH_k^m\bigl(\Sphere{N - 1}\bigr) \otimes
 \setC \frac{1}{\sqrt{2\pi}} r^{-\frac{2 \dindex{k} + N - 2}{2} + {\rm i}\sigma} \, {\rm d}\sigma
 \]
 and may write Theorem~\ref{thm:spectral-decomposition} as
 \begin{gather*}
 2E_x + 2 \dindex{k} + N - 2
 = \sumoplus_{m \in \setN} \intoplus_\setR 2{\rm i}\sigma \, {\rm d}\sigma ,\qquad
 {\rm i}\, \id
 = \sumoplus_{m \in \setN} \intoplus_\setR{\rm i} \, {\rm d}\sigma, \\
 {\rm i}\enorm{x}^2 \Laplacian_k
 = \sumoplus_{m \in \setN} \intoplus_\setR \Paren*{-{\rm i} \Paren*{\sigma^2 + \Paren*{m + \frac{2 \dindex{k} + N - 2}{2}}^2}} \, {\rm d}\sigma.
 \end{gather*}
\end{Remark}

\begin{Corollary}\label{thm:exponential}
 Let $ k $ be a non-negative multiplicity function.
 The $($possibly unbounded$)$ normal operator
 \[
 \exp\Paren*{\frac{z_1}{{\rm i}} (2E_x + 2 \dindex{k} + N - 2) + z_2 + z_3 \enorm{x}^2 \Laplacian_k}
 \]
 on $ L^2\bigl(\setR^N, w_{k,0}(x) {\rm d}x\bigr) $ is well-defined
 for $ z_1 , z_2 , z_3 \in \setC $.
 In particular, the action of the dif\-ferential-difference operators in
 $ \mfrakg_{k,0}$ lifts to a unique unitary representation of $ \setR^3 $ on
 $ L^2\bigl(\setR^N,\allowbreak w_{k,0}(x) {\rm d}x\bigr) $, which is given by
 \[
 (t_1, t_2, t_3) \mapsto
 \exp\bigl(t_1 (2E_x + 2 \dindex{k} + N - 2) + {\rm i}t_2 +{\rm i}t_3 \enorm{x}^2 \Laplacian_k\bigr).
 \]
\end{Corollary}

\begin{proof}
 Since $ \mfrakg_{k,0}$ admits joint spectral decomposition
 (see Theorem~\ref{thm:spectral-decomposition}), the (possibly unbounded) normal operator
 \[
 \phi\Paren*{\frac{1}{\rm i} (2E_x + 2 \dindex{k} + N - 2), 1, \enorm{x}^2 \Laplacian_k}
 \]
 on $ L^2\bigl(\setR^N, w_{k,0}(x) {\rm d}x\bigr) $ is defined for any Borel measurable
 function $ \map{\phi}{\setC^3}{\setC} $ by means of the functional calculus.
 The former assertion is shown by setting $ \phi(w_1, w_2, w_3)
 = \exp(z_1 w_1 + z_2 w_2 + z_3 w_3) $.
 The latter assertion is a consequence of Stone's theorem.
\end{proof}

For the operators in Theorem~\ref{thm:exponential}, the part involving $ z_2 $
contributes only as a scalar multiple of the identity.
The subsequent two subsections are devoted to the analysis of the parts
involving~$ z_1 $ and $ z_3 $.

\subsection[The unitary group with infinitesimal generator 2Ex + 2⟨k⟩ + N - 2]{The unitary group with infinitesimal generator $\boldsymbol{2E_x + 2 \dindex{k} + N - 2}$}

In this subsection, we consider the unitary group with infinitesimal generator
$ 2E_x + 2 \dindex{k} + N - 2 $, regarded as a skew-adjoint operator
on $ L^2\bigl(\setR^N, w_{k,0}(x) {\rm d}x\bigr) $ based on Theorem~\ref{thm:spectral-decomposition}.

\begin{Proposition}\label{thm:boundedness-1}
 Let $ k $ be a non-negative multiplicity function.
 For $ z \in \setC $, the $($possibly unbounded$)$ normal operator
 \[
 \exp\Paren*{\frac{z}{\rm i} (2E_x + 2 \dindex{k} + N - 2)}
 \]
 on $ L^2\bigl(\setR^N, w_{k,0}(x) {\rm d}x\bigr) $ is bounded if and only if
 $ \RePart z = 0 $, and in this case, this operator is unitary.
\end{Proposition}

\begin{proof}
 By Theorem~\ref{thm:spectral-decomposition}, the ope\-ra\-tor in the assertion
 corresponds to the multiplication operator \smash{$ \id_{L^2(\Sphere{N - 1}, w_k(\omega) {\rm d}\omega)} \otimeshat {\rm e}^{2z\sigma} $} on
 $ L^2\bigl(\Sphere{N - 1}, w_k(\omega) {\rm d}\omega\bigr) \otimeshat L^2(\setR, {\rm d}\sigma) $
 via the unitary operator \linebreak \smash{$ \id_{L^2(\Sphere{N - 1}, w_k(\omega) {\rm d}\omega)} \otimeshat (\mscrF \circ U_{N, k}) $}.
 Now the assertion follows since the function $ \sigma \mapsto {\rm e}^{2z\sigma} $
 on $ \setR $ is bounded if and only if $ \RePart z = 0 $, and in this case,
 $ \big|{\rm e}^{2z\sigma} \big| = 1 $.
\end{proof}

\begin{Theorem}\label{thm:scaling}
 Let $ k $ be a non-negative multiplicity function.
 For $ z = {\rm i}t $ with $ t \in \setR $, the unitary operator on
 $ L^2\bigl(\setR^N, w_{k,0}(x) {\rm d}x\bigr) $ in Theorem~{\rm\ref{thm:boundedness-1}}
 is given by
 \[
 \exp(t (2E_x + 2 \dindex{k} + N - 2)) F(x)
 = {\rm e}^{(2 \dindex{k} + N - 2)t} F\bigl({\rm e}^{2t} x\bigr).
 \]
\end{Theorem}

\begin{proof}
 We continue our discussion following the proof of Theorem~\ref{thm:boundedness-1}.
 The multiplication operator $ {\rm e}^{2{\rm i}t\sigma} $ on $ L^2(\setR, {\rm d}\sigma) $
 corresponds to the translation operator $ g \mapsto g((\blank) + 2t) $
 on $ L^2(\setR,{\rm d}s) $ via~$ \mscrF^{-1} $, which in turn corresponds to
 the scaling operator $ f \mapsto {\rm e}^{(2 \dindex{k} + N - 2)t} f\bigl({\rm e}^{2t} (\blank)\bigr) $
 on ${ L^2\bigl(\setRp, r^{2 \dindex{k} + N - 3} {\rm d}r\bigr) }$ via~$ U_{N, k}^{-1} $.
 Hence, the assertion holds.
\end{proof}

\subsection[The operator semigroup with infinitesimal generator |x|\^{}2 Delta\_k]{The operator semigroup with infinitesimal generator $\boldsymbol{\enorm{x}^2\Laplacian_k}$}

In this subsection, we consider the operator semigroup with infinitesimal generator
$ \enorm{x}^2 \Laplacian_k $, regarded as a self-adjoint operator on
$ L^2\bigl(\setR^N, w_{k,0}(x) {\rm d}x\bigr) $ based on Theorem~\ref{thm:spectral-decomposition}.

\begin{Proposition}\label{thm:boundedness-2}
 Let $ k $ be a non-negative multiplicity function.
 For $ z \in \setC $, the $($possibly~un\-bounded$)$ normal operator
$
 \exp\bigl(z \enorm{x}^2 \Laplacian_k\bigr)
$
 on $ L^2\bigl(\setR^N, w_{k,0}(x) {\rm d}x\bigr) $ is bounded if and only if ${\RePart z \geq 0}$, and unitary if and only if $ \RePart z = 0 $.
\end{Proposition}

\begin{proof}
 By Theorem~\ref{thm:spectral-decomposition}, the operator in the assertion
 corresponds to the multiplication operator
 \[
 \sideset{}{^\oplus}\sum_{m \in \mathbb{N}}
 \id_{\mcalH_k^m(\Sphere{N - 1})}
 \otimes {\rm e}^{-z (\sigma^2 + (m + \frac{2 \dindex{k} + N - 2}{2})^2)}\] on
 $ L^2\bigl(\Sphere{N - 1}, w_k(\omega) {\rm d}\omega\bigr) \otimeshat L^2(\setR, {\rm d}\sigma) $
 via the unitary operator
 $ \id_{L^2(\Sphere{N - 1}, w_k(\omega) {\rm d}\omega)} \otimeshat (\mscrF \circ U_{N, k}) $.
 Now the assertion follows since the function
\[
\sigma \mapsto {\rm e}^{-z (\sigma^2 + (m + \frac{2 \dindex{k} + N - 2}{2})^2)}
\]
 on $ \setR $ is bounded if and only if $ \RePart z \geq 0 $, and has modulus one
 if and only if $ \RePart z = 0 $.
\end{proof}

We consider the integral kernel formula for the operator semigroup
\smash{$ \bigl(\exp\bigl(z \enorm{x}^2 \Laplacian_k\bigr)\bigr)_{\RePart z \geq 0} $}.
For~this purpose, we first focus on the radial part
$ (E_r - m)(E_r + m + 2 \dindex{k} + N - 2) $ of~$ \enorm{x}^2 \Laplacian_k $~(see Theorem~\ref{thm:radial-parts-zero}).

\begin{Fact}[{\cite[Section~IX.7, Example~3]{MR493420}}]\label{thm:heat-kernel}
 The operator semigroup \smash{$ \bigl(\exp\bigl(z \bigl(\frac{\rm d}{{\rm d}s}\bigr)^2\bigr)\bigr)_{\RePart z \geq 0} $}
 on $ L^2(\setR,{\rm d}s) $ admits the integral kernel formula
 \begin{equation}
 \exp\Paren*{z \Paren*{\frac{\rm d}{{\rm d}s}}^2} g(s)
 = \frac{1}{\sqrt{4\pi z}}
 \int_\setR \exp\Paren*{-\frac{(s - s')^2}{4z}} g(s') \, {\rm d}s'
 \label{eq:heat-kernel}
 \end{equation}
 in the following sense. Here, we take the branch of $ \sqrt{z} $ such that
 $ \sqrt{z} > 0 $ when $ z > 0 $.
 \begin{enumarabicp}
 \item For $ z \in \setC $ with $ \RePart z > 0 $ and $ g \in L^2(\setR,{\rm d}s') $,
 the integrand in the right-hand side of \eqref{eq:heat-kernel} is
 integrable for all $ s \in \setR $, and this integral as a function
 of $ s $ gives \smash{$ \exp\bigl(z \bigl(\frac{\rm d}{{\rm d}s}\bigr)^2\bigr) g $}.
 \item For $ z \in \setC $ with $ \RePart z = 0 $ and $ z \neq 0 $ and
 $ g \in \bigl(L^1 \cap L^2\bigr)(\setR,{\rm d}s') $, the integrand in the right-hand
 side of \eqref{eq:heat-kernel} is integrable for all $ s \in \setR $,
 and this integral as a function of $ s $ gives~\smash{$ \exp\bigl(z \bigl(\frac{\rm d}{{\rm d}s}\bigr)^2\bigr) g $}.
 \end{enumarabicp}
\end{Fact}

\begin{Theorem}\label{thm:integral-kernel-rad}
 Let $ k $ be a non-negative multiplicity function and $ m \in \setN $.
 The operator semigroup
 $ \faml{\exp(z (E_r - m)(E_r + m + 2 \dindex{k} + N - 2))}{\RePart z \geq 0} $
 on $ L^2\bigl(\setRp, r^{2 \dindex{k} + N - 3} {\rm d}r\bigr) $ admits the integral kernel formula
 \begin{gather}
 \exp(z (E_r - m)(E_r + m + 2 \dindex{k} + N - 2)) f(r)
 = \int_{\setRp}\!\!\! K_k^{(m)}\bigl(r, r'; z\bigr) f\bigl(r'\bigr) r'{}^{2 \dindex{k} + N - 3} \, {\rm d}r',\!\!\!\!
 \label{eq:integral-kernel-of-radial-parts}
 \end{gather}
 where
 \begin{gather*}
 K_k^{(m)}(r, r'; z) \\
\qquad = \frac{1}{\sqrt{4\pi z}}
 \exp\Paren*{-z \Paren*{m + \frac{2 \dindex{k} + N - 2}{2}}^2}
 \exp\Paren*{-\frac{(\log r - \log r')^2}{4z}} \bigl(rr'\bigr)^{-\frac{2 \dindex{k} + N - 2}{2}}
 \end{gather*}
 for $ r , r' \in \setRp $, in the following sense. Here, we take the
 branch of $ \sqrt{z} $ such that $ \sqrt{z} > 0 $ when~${ z > 0} $.
 \begin{enumarabicp}
 \item For $ z \in \setC $ with $ \RePart z > 0 $ and
 $ f \in L^2\bigl(\setRp, r'{}^{2 \dindex{k} + N - 3} {\rm d}r'\bigr) $,
 the integrand in the right-hand side of \eqref{eq:integral-kernel-of-radial-parts}
 is integrable for all $ r \in \setRp $, and this integral as a function
 of $ r $ gives~$ {\exp(z (E_r - m)(E_r + m + 2 \dindex{k} + N - 2)) f} $.
 \item For $ z \in \setC $ with $ \RePart z = 0 $ and $ z \neq 0 $ and
 $ f \in \bigl(L^1 \cap L^2\bigr)\bigl(\setRp, r'{}^{2 \dindex{k} + N - 3} {\rm d}r'\bigr) $,
 the integrand in the right-hand side of \eqref{eq:integral-kernel-of-radial-parts}
 is integrable for all $ r \in \setRp $, and this integral as a function
 of $ r $ gives $ \exp(z (E_r - m)(E_r + m + 2 \dindex{k} + N - 2)) f $.
 \end{enumarabicp}
\end{Theorem}

\begin{proof}
 As follows from the proof of Theorem~\ref{thm:spectral-decomposition-rad},
 $ (E_r - m)(E_r + m + 2 \dindex{k} + N - 2) $ corresponds to
 \smash{$ \bigl(\frac{\rm d}{{\rm d}s}\bigr)^2 - \bigl(m + \frac{2 \dindex{k} + N - 2}{2}\bigr)^2 $} via the unitary operator
 $ \map{U_{N, k}}{L^2(\setRp, r^{2 \dindex{k} + N - 3} {\rm d}r)}{L^2(\setR,{\rm d}s)} $.
 Now the assertion follows from Theorem~\ref{thm:heat-kernel}.
\end{proof}

We then combine this result with the spherical part. Recall that \smash{$ P_k^{(m)} $}
denotes the Poisson kernel of the space of $ k $-spherical harmonics of
degree $ m $ (see \eqref{eq:poisson-kernel} for the definition), and that
$ C_m^\nu $ \big(resp.\ $ \widetilde{C}_m^\nu $\big) denotes the Gegenbauer polynomial
(resp.\ the renormalized Gegenbauer polynomial) (see \eqref{eq:gegenbauer} and
\eqref{eq:renormalized-gegenbauer} for the definitions).

\begin{Lemma}\label{thm:uniform-norm-of-poisson-kernel}
 Fix a non-negative multiplicity function $ k $. Then, the uniform norm
 of the Poisson kernel \smash{$ P_k^{(m)} $} satisfies
 \[
 \sup_{\omega, \omega' \in \Sphere{N - 1}} \big|P_k^{(m)}(\omega, \omega')\big|
 = O\bigl(m^{2 \dindex{k} + N - 2}\bigr),
 \qquad m \to \infty.
 \]
\end{Lemma}

\begin{proof}
 In the case $ N = 1 $, we have $ P_k^{(m)} = 0 $ for $ m \geq 2 $
 (see Theorem~\ref{rem:poisson-kernel-of-1-dim}), so that the assertion holds trivially.
 We now consider the case $ N \geq 2 $.
 By Theorems~\ref{thm:poisson-kernel} and \ref{thm:dunkls-intertwining-operator-as-integral}, we have
 \begin{align*}
 P_k^{(m)}(\omega, \omega')
 &= V_k\bigl(\widetilde{C}_m^{\frac{2 \dindex{k} + N - 2}{2}}(\innprod{\blank}{\omega'})\bigr)(\omega)
 = \int_{\setR^N} \widetilde{C}_m^{\frac{2 \dindex{k} + N - 2}{2}}(\innprod{\xi}{\omega'})
 \, {\rm d}\mu_{k, \omega}(\xi),
 \end{align*}
 where $ \mu_{k, x} $ is a probability Borel measure on $ \setR^N $
 whose support is contained in the unit ball
 $ \bigl\{\xi \in \setR^N\mid \enorm{\xi} \leq 1\bigr\} $.
 Hence, we have
 \begin{align*}
\big|P_k^{(m)}(\omega, \omega')\big|
 &\leq \int_{\setR^N} \big|\widetilde{C}_m^{\frac{2 \dindex{k} + N - 2}{2}}(\innprod{\xi}{\omega'})\big|
 \, {\rm d}\mu_{k, \omega}(\xi) = \sup_{t \in [-1, 1]} \big|\widetilde{C}_m^{\frac{2 \dindex{k} + N - 2}{2}}(t)\big|.
 \end{align*}
 Note that $ \frac{2 \dindex{k} + N - 2}{2} \geq 0 $ since $ N \geq 2 $.
 For $ \nu \in \setRp $, it is known (see \cite[p.~302]{MR1688958}) that
 \[
 \sup_{t \in [-1, 1]} \abs{C_m^\nu(t)}
 = C_m^\nu(1)
 = \frac{\Gamma(m + 2\nu)}{m! \Gamma(2\nu)},
 \]
 which implies
 \[
 \sup_{t \in [-1, 1]} \big|\widetilde{C}_m^\nu(t)\big|
 = \frac{m + \nu}{\nu} \frac{\Gamma(m + 2\nu)}{m! \Gamma(2\nu)}
 = O\bigl(m^{2\nu}\bigr).
 \]
 This also holds for $ \nu = 0 $ since $ \widetilde{C}_m^0 $ is defined
 by the limit formula \eqref{eq:limit-formula}.
 The assertion follows by applying this estimate to the case
 \smash{$ \nu = \frac{2 \dindex{k} + N - 2}{2} $}.
\end{proof}

Recall from \eqref{eq:volume} that $ \vol_k\bigl(\Sphere{N - 1}\bigr) $ denotes the volume
of the sphere $ \Sphere{N - 1} $ with respect to the measure
$ w_k(\omega) {\rm d}\omega $.

\begin{Theorem}\label{thm:integral-kernel}
 Let $ k $ be a non-negative multiplicity function.
 The operator semigroup \[\bigl(\exp\bigl(z \enorm{x}^2 \Laplacian_k\bigr)\bigr)_{\RePart z \geq 0} \]
 on $ L^2\bigl(\setR^N, w_{k,0}(x) {\rm d}x\bigr) $ admits the integral kernel formula
 \begin{equation}
 \exp\bigl(z \enorm{x}^2 \Laplacian_k\bigr) F(x)
 = \int_{\setR^N} K_k\bigl(x, x'; z\bigr) F\bigl(x'\bigr) w_{k,0}\bigl(x'\bigr) \, {\rm d}x',
 \label{eq:integral-kernel}
 \end{equation}
 where
 \begin{align*}
 K_k(r\omega, r'\omega'; z)
 ={}&\frac{1}{\sqrt{4\pi z}}
 \exp\Paren*{-\frac{(\log r - \log r')^2}{4z}}
 \bigl(rr'\bigr)^{-\frac{2 \dindex{k} + N - 2}{2}} \\
 &
 \times \frac{1}{\vol_k\bigl(\Sphere{N - 1}\bigr)} \sum_{m = 0}^{\infty}
 \exp\Paren*{-z \Paren*{m + \frac{2 \dindex{k} + N - 2}{2}}^2}
 P_k^{(m)}(\omega, \omega')
 \end{align*}
 for $ r, r' \in \setRp $ and $ \omega, \omega' \in \Sphere{N - 1} $,
 in the following sense. Here, we take the branch of $ \sqrt{z} $
 such that $ \sqrt{z} > 0 $ when $ z > 0 $.

 For $ z \in \setC $ with $ \RePart z > 0 $ and
 $ F \in L^2\big(\setR^N, w_{k,0}\bigl(x'\bigr) {\rm d}x'\big) $, the integrand in the right-hand
 side of~\eqref{eq:integral-kernel} is integrable for all
 $ x \hspace{-0.41pt}\in\hspace{-0.41pt} \setR^N \setminus \setenum{0} $, and this integral as a function of
 $ x $ gives~${ \exp\bigl(z \enorm{x}^2 \Laplacian_k\bigr) F} $.
\end{Theorem}

\begin{Remark}
 When $ k = 0 $, we have
 \[
 \vol_0\bigl(\Sphere{N - 1}\bigr) = \vol\bigl(\Sphere{N - 1}\bigr)
 = \frac{2\pi^{\frac{N}{2}}}{\Gamma\bigl(\frac{N}{2}\bigr)}, \qquad \text{with} \quad P_0^{(m)}(\omega, \omega')
 = \widetilde{C}_m^{\frac{N - 2}{2}}(\innprod{\omega}{\omega'})
 \]
 (see Theorem~\ref{rem:classical-poisson-kernel}), so that
 the integral kernel formula in Theorem~\ref{thm:integral-kernel} reduces to
 \begin{align*}
 K_0(r\omega, r'\omega'; z)
 &= \frac{1}{\sqrt{4\pi z}}
 \exp\Paren*{-\frac{\bigl(\log r - \log r'\bigr)^2}{4z}}
 \bigl(rr'\bigr)^{-\frac{N - 2}{2}} \\
 &\qquad
 \times \frac{\Gamma\bigl(\frac{N}{2}\bigr)}{2\pi^{\frac{N}{2}}} \sum_{m = 0}^{\infty}
 \exp\Paren*{-z \Paren*{m + \frac{N - 2}{2}}^2}
 \widetilde{C}_m^{\frac{N - 2}{2}}(\innprod{\omega}{\omega'}).
 \end{align*}
\end{Remark}

\begin{proof}[Proof of Theorem~\ref{thm:integral-kernel}]
 Fix $ r \in \setRp $ and $ \omega \in \Sphere{N - 1} $. Then, the function
 \[
 r' \mapsto
 \exp\Paren*{-\frac{(\log r - \log r')^2}{4z}} \bigl(rr'\bigr)^{-\frac{2 \dindex{k} + N - 2}{2}}
 \]
 is square-integrable with respect to the measure $ r'{}^{2 \dindex{k} + N - 3} {\rm d}r' $
 and the infinite series
 \[
 \sum_{m = 0}^{\infty}
 \exp\Paren*{-z \Paren*{m + \frac{2 \dindex{k} + N - 2}{2}}^2}
 P_k^{(m)}(\omega, \omega')
 \]
 absolutely converges with respect to the uniform norm for
 $ \omega' \in \Sphere{N - 1} $ by Theorem~\ref{thm:uniform-norm-of-poisson-kernel},
 so that the equation
 \begin{align*}
 K_k(r\omega, r'\omega'; z)
 ={}& \frac{1}{\sqrt{4\pi z}}
 \exp\Paren*{-\frac{(\log r - \log r')^2}{4z}}
 \bigl(rr'\bigr)^{-\frac{2 \dindex{k} + N - 2}{2}} \\
 &
 \times \frac{1}{\vol_k\bigl(\Sphere{N - 1}\bigr)} \sum_{m = 0}^{\infty}
 \exp\Paren*{-z \Paren*{m + \frac{2 \dindex{k} + N - 2}{2}}^2}
 P_k^{(m)}(\omega, \omega') \\
 ={}& \frac{1}{\vol_k\bigl(\Sphere{N - 1}\bigr)}
 \sum_{m = 0}^{\infty} P_k^{(m)}(\omega, \omega') K_k^{(m)}(r, r'; z)
 \end{align*}
 (here, $ K_k^{(m)}(r, r'; z) $ is as defined in Theorem~\ref{thm:integral-kernel-rad})
 holds with respect to the topology of
 $ L^2\bigl(\setRp \times \Sphere{N - 1}, r'{}^{2 \dindex{k} + N - 3} w_k(\omega') {\rm d}r' {\rm d}\omega'\bigr)
 \cong L^2\bigl(\setR^N, w_{k,0}\bigl(x'\bigr) {\rm d}x'\bigr) $.

 By the result of the previous paragraph,
 for $ F \in L^2\bigl(\setR^N, w_{k,0}\bigl(x'\bigr) {\rm d}x'\bigr) $
 and $ x = r\omega \in \setR^N \setminus \setenum{0} $,
 the function $ x' \mapsto K_k\bigl(x, x'; z\bigr) F\bigl(x'\bigr) $ is integrable with respect to
 the measure $ w_{k,0}\bigl(x'\bigr) {\rm d}x' $ and
 \begin{align*}
 \int_{\setR^N} K_k\bigl(x, x'; z\bigr) F\bigl(x'\bigr) w_{k,0}\bigl(x'\bigr) \, {\rm d}x'
 ={}& \frac{1}{\vol_k\bigl(\Sphere{N - 1}\bigr)}
 \sum_{m = 0}^{\infty} \int_{\Sphere{N - 1}} \int_{\setRp}
 P_k^{(m)}\bigl(\omega, \omega'\bigr)
\\
 &\times K_k^{(m)}(r, r'; z) F\bigl(r'\omega'\bigr)
 r'{}^{2 \dindex{k} + N - 3} w_k\bigl(\omega'\bigr) \, {\rm d}r' {\rm d}\omega'.
 \end{align*}
 If $ F = p \otimes f $ with $ p \in \mcalH_k^l\bigl(\Sphere{N - 1}\bigr) $ and
 $ f \in L^2\bigl(\setRp, r'{}^{2 \dindex{k} + N - 3} {\rm d}r'\bigr) $,
 by \eqref{eq:poisson-kernel} and Theorem~\ref{thm:integral-kernel-rad}\,(1),
 we have
 \begin{gather*}
 \int_{\setR^N} K_k\bigl(x, x'; z\bigr) (p \otimes f)\bigl(x'\bigr) w_{k,0}\bigl(x'\bigr) \, {\rm d}x' \\
 \qquad= \frac{1}{\vol_k\bigl(\Sphere{N - 1}\bigr)}
 \sum_{m = 0}^{\infty} \int_{\Sphere{N - 1}} \int_{\setRp}
 P_k^{(m)}(\omega, \omega')
 K_k^{(m)}(r, r'; z) p(\omega') f(r')\\
 \phantom{\qquad=}{}\times
 r'{}^{2 \dindex{k} + N - 3} w_k(\omega') \, {\rm d}r' {\rm d}\omega' \\
 \qquad= \sum_{m = 0}^{\infty}
 \Paren*{
 \frac{1}{\vol_k\bigl(\Sphere{N - 1}\bigr)}
 \int_{\Sphere{N - 1}} P_k^{(m)}(\omega, \omega') p(\omega') w_k(\omega') \, {\rm d}\omega'
 }\\
 \phantom{\qquad=}{}\times
 \Paren*{\int_{\setRp} K_k^{(m)}(r, r'; z) f(r') r'{}^{2 \dindex{k} + N - 3} \, {\rm d}r'} \\
 \qquad= p(\omega) \Paren*{\int_{\setRp} K_k^{(l)}(r, r'; z) f(r') r'{}^{2 \dindex{k} + N - 3} \, {\rm d}r'} \\
 \qquad= p(\omega) \exp(z (E_r - l)(E_r + l + 2 \dindex{k} + N - 2)) f(r) \\
\qquad= \exp\bigl(z \enorm{x}^2 \Laplacian_k\bigr) (p \otimes f)(x).
 \end{gather*}
 Hence, \eqref{eq:integral-kernel} holds in this case.

 Let $ F \in L^2\bigl(\setR^N, w_{k,0}\bigl(x'\bigr) {\rm d}x'\bigr) $ and take a sequence
 $ \faml{F_j}{j \in \setN} $ in
 \[ \mcalH\bigl(\Sphere{N - 1}\bigr) \otimes L^2\bigl(\setRp, r'{}^{2 \dindex{k} + N - 3} {\rm d}r'\bigr) \]
 such that $ F_j \to F $ in $ L^2\bigl(\setR^N, w_{k,0}\bigl(x'\bigr) {\rm d}x'\bigr) $.
 Then, since $ \exp\bigl(z \enorm{x}^2 \Laplacian_k\bigr) $ is a bounded operator on~${ L^2\bigl(\setR^N, w_{k,0}(x) {\rm d}x\bigr) }$ (see Theorem~\ref{thm:boundedness-2}),
 we have
 \[
 \exp\bigl(z \enorm{x}^2 \Laplacian_k\bigr) F_j
 \to \exp\bigl(z \enorm{x}^2 \Laplacian_k\bigr) F
 \qquad \text{in}\quad L^2\bigl(\setR^N, w_{k,0}(x) {\rm d}x\bigr) .
 \]
 On the other hand, for each $ x \in \setR^N \setminus \setenum{0} $,
 the function $ x' \mapsto K_k\bigl(x, x'; z\bigr) $ is square-integrable with respect to
 the measure $ w_{k,0}\bigl(x'\bigr) {\rm d}x' $, so we have
 \[
 \int_{\setR^N} K_k\bigl(x, x'; z\bigr) F_j\bigl(x'\bigr) w_{k,0}\bigl(x'\bigr) \, {\rm d}x'
 \to \int_{\setR^N} K_k\bigl(x, x'; z\bigr) F\bigl(x'\bigr) w_{k,0}\bigl(x'\bigr) \, {\rm d}x'.
 \]
 By the result of the previous paragraph, \eqref{eq:integral-kernel} holds for
 each $ F_j $. By taking the limit as $ j \to \infty $, we conclude that
 \eqref{eq:integral-kernel} also holds for $ F $.
\end{proof}

\begin{Remark}\label{rem:laguerre-semigroup}
 We recall the definition of the $ (k, a) $-generalized Laguerre semigroup
 \[ \faml{\mscrI_{k,a}(z)}{\RePart z \geq 0} \] from \cite[equation~(1.3)]{MR2956043}:
 \[
 \mscrI_{k,a}(z)
 = \exp\left(\frac{z}{{\rm i}} \bigl(\DiffEm{k, a} - \DiffEp{k, a}\bigr)\right)
 = \exp\left(\frac{z}{a} \bigl(\enorm{x}^{2 - a} \Laplacian_k - \enorm{x}^a\bigr)\right),
 \]
 and the definition of the $ (k, a) $-generalized Fourier transform
 $ \mscrF_{k,a} $ from \cite[equation~(5.2)]{MR2956043}:
 \begin{align*}
 \mscrF_{k,a}
 &= {\rm e}^{\frac{{\rm i}\pi}{2} \frac{2 \dindex{k} + a + N - 2}{a}}
 \mscrI_{k,a}\Paren*{\frac{{\rm i}\pi}{2}} = {\rm e}^{\frac{{\rm i}\pi}{2} \frac{2 \dindex{k} + a + N - 2}{a}}
 \exp\Paren*{\frac{\pi}{2} \bigl(\DiffEm{k, a} - \DiffEp{k, a}\bigr)} \\
 &= {\rm e}^{\frac{{\rm i}\pi}{2} \frac{2 \dindex{k} + a + N - 2}{a}}
 \exp\Paren*{\frac{{\rm i}\pi}{2a} \bigl(\enorm{x}^{2 - a} \Laplacian_k - \enorm{x}^a\bigr)}.
 \end{align*}

 These are not well-defined for $ a = 0 $. However, considering
 the ``renormalized'' $ (k, a) $-generalized Laguerre semigroup
$
 \mscrI_{k,a}(az)
 = \exp\bigl(z \bigl(\enorm{x}^{2 - a} \Laplacian_k - \enorm{x}^a\bigr)\bigr)
$
 and putting $ a = 0 $, we get the operator
$
 \exp\bigl(z \bigl(\enorm{x}^2 \Laplacian_k - 1\bigr)\bigr)
 = {\rm e}^{-z} \exp\bigl(z \enorm{x}^2 \Laplacian_k\bigr)$.
 By Theorem~\ref{thm:integral-kernel}, for $ z \in \setC $ with~${ \RePart z > 0 }$,
 the integral kernel of this operator is the function
 $ (x, x') \mapsto {\rm e}^{-z} K_k\bigl(x, x'; z\bigr) $.
\end{Remark}

\section[Closed-form expressions for the integral kernels in low-dimensional cases]{Closed-form expressions for the integral kernels\\ in low-dimensional cases}
\label{sec:closed-form}

In this section, we give a closed-form expression for the integral kernel
$ (x, x') \mapsto K_k\bigl(x, x'; z\bigr) $ of the operator $ \exp\bigl(z \enorm{x}^2 \Laplacian_k\bigr) $
($ \RePart z > 0 $), obtained in Theorem~\ref{thm:integral-kernel},
in the low-dimensional cases~${ N = 1 , 2 }$ and $ 4 $.
In the cases $ N = 2 $ and $ 4 $, we assume that $ k = 0 $, and show that
the integral kernel can be expressed in terms of the theta function.

\begin{Proposition}\label{thm:integral-kernel-of-1-dim}
 We consider the case $ N = 1 $. The reduced root system $ \mscrR $
 is taken to be~$ {\setenum{\alpha, -\alpha} }$ with $ \alpha \in \setRp $,
 and the non-negative multiplicity function $ k $ is identified with
 $ k_\alpha = k_{-\alpha} \in \setRzp $.
 In this case, for $ z \in \setC $ with $ \RePart z > 0 $, we have
 \begin{align*}
 K_k\bigl(x, x'; z\bigr)
 ={}& \frac{1}{2\alpha^{2k} \sqrt{4\pi z}}
 \exp\Paren*{-\frac{(\log \abs{x} - \log \abs{x'})^2}{4z}} \abs{xx'}^{-k + \frac{1}{2}}\\
 &\times
 \bigl({\rm e}^{-(k - \frac{1}{2})^2 z} + {\rm e}^{-(k + \frac{1}{2})^2 z} \sgn(xx')\bigr)
 \end{align*}
 for $ x , x' \in \setR \setminus \setenum{0} $.
 In particular, when $ k = 0 $, we have
 \[
 K_0\bigl(x, x'; z\bigr)
 = \frac{{\rm e}^{-\frac{z}{4}}}{2 \sqrt{4\pi z}}
 \exp\Paren*{-\frac{(\log \abs{x} - \log \abs{x'})^2}{4z}} \abs{xx'}^{\frac{1}{2}}
 (1 + \sgn(xx'))
 \]
 for $ x, x' \in \setR \setminus \setenum{0} $.
\end{Proposition}

\begin{proof}
 It follows from $ \vol_k\bigl(\Sphere{0}\bigr) = 2\alpha^{2k} $ and
 Theorem~\ref{rem:poisson-kernel-of-1-dim}.
\end{proof}

In the following, we consider only the case $ k = 0 $.
We recall the definition of the theta function:
\[
 \vartheta(v, \tau)
 = \sum_{m = -\infty}^{\infty} \exp\bigl({\rm i}\pi\tau m^2 + 2{\rm i}\pi mv\bigr)
 = 1 + 2 \sum_{m = 1}^{\infty} \exp\bigl({\rm i}\pi\tau m^2\bigr) \cos 2\pi mv.
\]

\begin{Proposition}\label{thm:integral-kernel-of-2-dim}
 In the case $ N = 2 $, for $ z \in \setC $ with $ \RePart z > 0 $, we have
 \[
 K_0\bigl(r\omega, r'\omega'; z\bigr)
 = \frac{1}{2\pi \sqrt{4\pi z}}
 \exp\Paren*{-\frac{(\log r - \log r')^2}{4z}}
 \vartheta\Paren*{\frac{1}{2\pi} \arccos \innprod{\omega}{\omega'}, \frac{\rm i}{\pi} z}
 \]
 for $ r , r' \in \setRp $ and $ \omega , \omega' \in \Sphere{1} $.
 Equivalently, we have
 \[
 K_0\bigl(r{\rm e}^{{\rm i}\phi}, r'{\rm e}^{{\rm i}\phi'}; z\bigr)
 = \frac{1}{2\pi \sqrt{4\pi z}}
 \exp\Paren*{-\frac{\bigl(\log r - \log r'\bigr)^2}{4z}}
 \vartheta\Paren*{\frac{1}{2\pi} \bigl(\phi - \phi'\bigr), \frac{\rm i}{\pi} z}
 \]
 for $ r, r' \in \setRp $ and $ \phi, \phi' \in \setR $.
\end{Proposition}

\begin{proof}
 Since
 \[
 \widetilde{C}_m^0(t)
 = \begin{cases}
 1, & m = 0, \\
 2 T_m(t), &m \geq 1,
 \end{cases}
 \]
 where $ T_m $ denotes the Chebyshev polynomial of the first kind,
 which is characterized by the formula and $ T_m(\cos \theta) = \cos m\theta $,
 we have
 \begin{align*}
 \sum_{m = 0}^{\infty} \exp\bigl(-zm^2\bigr) \widetilde{C}_m^0(\cos \theta)
 &= 1 + 2 \sum_{m = 1}^{\infty} \exp\bigl(-zm^2\bigr) T_m(\cos \theta) \\
 &= 1 + 2 \sum_{m = 1}^{\infty} \exp\bigl(-zm^2\bigr) \cos m\theta = \vartheta\Paren*{\frac{1}{2\pi} \theta, \frac{\rm i}{\pi} z}.
 \end{align*}
 Hence, we have
 \begin{align*}
 K_0\bigl(r\omega, r'\omega'; z\bigr)
 &= \frac{1}{\sqrt{4\pi z}}
 \exp\Paren*{-\frac{\bigl(\log r - \log r'\bigr)^2}{4z}}
 \times \frac{1}{2\pi} \sum_{m = 0}^{\infty} \exp\bigl(-zm^2\bigr) \widetilde{C}_m^0\bigl(\innprod{\omega}{\omega'}\bigr) \\
 &= \frac{1}{2\pi \sqrt{4\pi z}}
 \exp\Paren*{-\frac{\bigl(\log r - \log r'\bigr)^2}{4z}}
 \vartheta\Paren*{\frac{1}{2\pi} \arccos \innprod{\omega}{\omega'}, \frac{\rm i}{\pi} z}.\tag*{\qed}
 \end{align*}\renewcommand{\qed}{}
\end{proof}

\begin{Proposition}\label{thm:integral-kernel-of-4-dim}
 In the case $ N = 4 $, for $ z \in \setC $ with $ \RePart z > 0 $, we have
 \begin{align*}
 K_0\bigl(r\omega, r'\omega'; z\bigr)
 ={}& -\frac{1}{8\pi^3 \sqrt{4\pi z}}
 \exp\Paren*{-\frac{(\log r - \log r')^2}{4z}} \bigl(rr'\bigr)^{-1}
\\ &\times
\bigl(1 - \innprod{\omega}{\omega'}^2\bigr)^{-\frac{1}{2}} \Pdiff{\vartheta}{v} \Paren*{\frac{1}{2\pi} \arccos \innprod{\omega}{\omega'}, \frac{\rm i}{\pi} z}
 \end{align*}
 for $ r, r' \in \setRp $ and $ \omega, \omega' \in \Sphere{3} $.
 Here, we take the branch of $ \arccos \innprod{\omega}{\omega'} $ such that
 $ \arccos \innprod{\omega}{\omega'} \in [0, \pi] $.
\end{Proposition}

\begin{proof}
 Since
$
 \widetilde{C}_m^1(t)
 = (m + 1) C_m^1(t)
 = (m + 1) U_m(t)$,
 where $ U_m $ denotes the Chebyshev polynomial of the second kind, which is
 characterized by the formula $ U_m(\cos \theta)
 = (\sin (m + 1)\theta)/(\sin \theta) $, we have
 \begin{align*}
 \begin{split}
 \sum_{m = 0}^{\infty} \exp\bigl(-z(m + 1)^2\bigr) \widetilde{C}_m^1(\cos \theta)
 &= \sum_{m = 0}^{\infty} \exp\bigl(-z(m + 1)^2\bigr) \cdot (m + 1) \frac{\sin (m + 1)\theta}{\sin \theta} \\
 &= \frac{1}{\sin \theta} \sum_{m = 1}^{\infty} \exp\bigl(-zm^2\bigr) \cdot m \sin m\theta \\
 &= -\frac{1}{4\pi \sin \theta}
 \Pdiff{\vartheta}{v} \Paren*{\frac{1}{2\pi} \theta, \frac{\rm i}{\pi} z}.
 \end{split}
 \end{align*}
 Hence, we have
 \begin{gather*}
 K_0(r\omega, r'\omega'; z) \\
\qquad = \frac{1}{\sqrt{4\pi z}}
 \exp\Paren*{-\frac{\bigl(\log r - \log r'\bigr)^2}{4z}} \bigl(rr'\bigr)^{-1}
 \times \frac{1}{2\pi^2} \sum_{m = 0}^{\infty} \exp(-z(m + 1)^2) \widetilde{C}_m^1\bigl(\innprod{\omega}{\omega'}\bigr) \\
\qquad = -\frac{1}{8\pi^3 \sqrt{4\pi z}}
 \exp\Paren*{-\frac{\bigl(\log r - \log r'\bigr)^2}{4z}} \bigl(rr'\bigr)^{-1}
\\
 \phantom{\qquad = }{}\times
\bigl(1 - \innprod{\omega}{\omega'}^2\bigr)^{-\frac{1}{2}} \Pdiff{\vartheta}{v} \Paren*{\frac{1}{2\pi} \arccos \innprod{\omega}{\omega'}, \frac{\rm i}{\pi} z}.\tag*{\qed}
 \end{gather*}\renewcommand{\qed}{}
\end{proof}

\subsection*{Acknowledgements}

The author presented this work at the workshop ``Intertwining operators and
geometry'' (January 20--24, 2025) held at the Institut Henri Poincar\'e, Paris.
He gratefully acknowledges the organizers Professor Jan Frahm and Professor
Angela Pasquale for the invitation and their kindness during the workshop.
The author would also like to express his gratitude to Professor Toshiyuki
Kobayashi for his guidance, encouragement, and insightful discussions.
Finally, the author thanks the anonymous referees for their helpful comments,
which have led him to improve this paper.

This work was supported by World-leading Innovative Graduate Study for Frontiers
of Mathematical Sciences and Physics (WINGS-FMSP), The University of Tokyo,
and JSPS KAKENHI Grant Number JP25KJ0914.

\pdfbookmark[1]{References}{ref}
\LastPageEnding

\end{document}